# Living conditions:
# classification of households using the Kohonen algorithm

Sophie Ponthieux[*] and Marie Cottrell[**]

[*] INSEE, Division "Conditions de vie des ménages" – Timbre F340
18 Bd A.Pinard – F 75675 Paris Cedex 14 – France
phone: + 33 1 41 17 38 93 / fax: + 33 1 41 17 63 17
E-mail: sophie.ponthieux@insee.fr

[**] SAMOS-MATISSE, Université de Paris 1
90 rue de Tolbiac – F 75634 Paris Cedex 13
phone-fax: + 33 1 44 07 89 22
E-mail: cottrell@univ-paris1.fr

# Living conditions:
# classification of households using the Kohonen algorithm

Sophie Ponthieux and Marie Cottrell

**Abstract:**

In the analysis of poverty and social exclusion, indicators of living conditions are some interesting non-monetary complements to the usual measurements in terms of current or annual income. Living conditions depend in fact on longer term factors than income, and provide further information on households' actual resources that allow to compare more accurately between living standards. But in counterpart, a difficulty comes from the qualitative nature of the information, and the large number of dimensions and items that may be taken into account; in other words, living conditions are difficult to "measure". A consequence is that very often, the information is either used only partly, or reduced into a global score of (bad) living conditions, that results from counting "negative" items, and the qualitative dimension is lost. In this paper, we propose to use the Kohonen algorithm first to describe how the elements of living conditions are combined, and secondly to classify households according to their living conditions. The main interest of a classification is to make appear not only quantitative differences in the "levels" of living conditions, but also qualitative differences within similar "levels".

# 1 Introduction

Since the 1970s in United-Kingdom, more recently in France, the use of non-monetary indicators in the analysis of poverty is developing, and poverty or exclusion are studied both in terms of income and in terms of living conditions (Townsend, 1979; Nolan & Whelan, 1996). Living conditions include a great number of domains. Dickes (1994) lists ten of them: dwelling, durables, food, clothing, financial resources, health, social relations, leisure, education and work. Not all the existing studies include this complete set of domains. The choice of including or not one of those domains may be based on two arguments: in the first one, the main hypothesis is that the subjects are rational in their behavior, which leads to select only the domains where privations are supposed to decrease or disappear when the financial resources increase (Mack & Lansley, 1984); the second one is based on the notion of "standard" (Townsend, 1989; Dickes, 1994), and leads to consider any domain as soon as all the subjects are - at least potentially – involved. Each domain is described by a list of items, that can be, for example, a characteristic of the dwelling ("*is there central heating/bath or shower/etc. in the house*"?), a particular consumer good ("*do you have a color TV/a car/a washing machine/etc.*"?), the ability to afford ("*new clothes/a week's holiday away/etc.*"?); for each item, the respondent for the household declares to have or have not, to be able to afford or not.

Whatever the option in choosing the domains, the main difficulty is to deal with a great quantity of information (depending on the number of dimensions retained and the number of items within each dimension) that is mainly qualitative. One possibility is to construct a "score" of (bad) living conditions, which corresponds to the total number of "negative" items (for example, Lollivier & Verger, 1997). The interest of this approach is that high scores will always characterize households who deal with cumulative difficulties; nevertheless it will not allow to distinguish between households who have smaller scores, but difficulties of completely different nature; and since there is no method for weighting the items, not having a TV set has exactly the same value as not being able to afford buying meat/chicken/fish every second day, even though it is clearly not the same nature of hardship.

An alternative to using a score is to classify households into groups that would be consistent both in terms of living conditions "level" and in terms of (bad) living conditions nature. For this purpose, several classification methods can be used; in what follows, we use a classification algorithm relating to self-organizing maps, the Kohonen algorithm. The data are from the third wave of the French part of the European community households panel.

In a first step, we classify only the modalities, in order to obtain a good description of living conditions characteristics: how they are combined, *i.e.* what are the most frequent associations between the items modalities. The result show some neat oppositions not only between "positive" and "negative" modalities but also between groups of items, that could suggest some kind of qualitative gradation in living conditions hardship. In a second step, we classify the observations, searching for a consistent grouping of households described only by their living conditions. The classifications obtained tend to confirm that beyond differences in the "level" of (bad) living conditions, there are significant within level differences in the nature of the difficulties.

## 2 Data, variables of the analysis and method

The data source is the French part of the European Community Households Panel, here in its third wave (year 1996). It provides the required detailed information about material living conditions (dwelling, environment, durable goods, deprivations); we also know the household's income, and whether the household's respondent considers that the monetary resources allow to live from "very comfortably" to "with great difficulty" ("subjective" living conditions after). Finally, the source covers three dimensions of households' poverty: in terms of monetary income, in terms of material living conditions, and in terms of "subjective" living conditions; this will allow to characterize the classes of living conditions from the two other points of view.

The observations are households. For the classifications, we use only their characteristics relating to material living conditions. Living conditions are described by 10 items about the dwelling (5 about convenience and 5 about problems), 4 items about environmental topics, 6 items for the durables and 6 items about deprivations, a total of 26 dummy variables, that is 52 modalities (detailed in Appendix, table A1). For each item, the "negative" modality (having a problem, not having an item, not being able to afford) is always coded "1" *vs.* "0" for the "neutral" modality.

The observations are also described by a set of general characteristics (Appendix, table A2): type of household, average age of the adults (persons aged 17 years and over), number of children under 17 years, current income per consumption unit, and type and location of the dwelling. In addition, we also characterize households by an indicator of monetary poverty, an indicator of subjective living conditions, and scores of material living conditions (total and partial – by domain -). These general characteristics are used only to compare between classes, not to construct them.

Only the observations with no missing variable for all these descriptors are kept for the analysis, that is 6458 households.

For the classifications, the information in input is not always under the same form: we use a response table to classify the characteristics and successively the partial scores (score by domain) then the coordinates (after a multiple correspondence analysis, MCA) to classify the observations. The classifications are done using the Kohonen algorithm (for an introduction to the algorithm and its applications to data analysis, see Kohonen 1984, 1993, 1995; Kaski, 1997; Cottrell & Rousset, 1997). The main property of the Kohonen algorithm is its property to preserve the topology of the data: after convergence, similar data are grouped into the same class or into neighbor classes. This feature allows to represent the proximity between data, as in a projection, in the Kohonen map. As a further treatment, the Kohonen classes can be clustered into a reduced number of macro-classes (which only contain neighbor Kohonen classes) by using a classical hierarchical classification.

**3. Classification of the modalities**

In this first step, the objective is to obtain a representation of the combination between the modalities of the whole set of items. We have tested successively a one-dimension and a two-dimensions classifications using a Kohonen algorithm inspired by the MCA (Cottrel, Letremy & Roy, 1994; Cottrell & Ibbou, 1995).

The results from the first one are represented along a ten classes string (Kohonen map of one dimension). We obtain seven classes of "negative" modalities, and three classes of "neutral" modalities. On the "negative" side, the first class (reading Figure 1 from the left) associates low standard dwelling (no bath or shower, no hot running water and no indoor toilet) and the absence of a telephone and TV set, which are very common items (owned by about 99% of the households). The next class groups together not having a car and not being able to afford one meal of meat/fish/chicken every second day, that is a very serious deprivation. The third class associates several deprivations and the absence of some "modern" durable goods (no micro-wave oven, no VCR - but it is to be related to the absence of a TV set -, and no dish-washer). The three next classes are characterized by the association of problems relating to the quality of he dwelling (dampness, shortage of space) and environmental problems (pollution, noise), the last class of this group adding the inability to afford one week's holiday away from home. The last of the "negative" classes contains only one characteristic which is the inability to replace broken or worn furniture. The interest of this classification is that it indicates a qualitative gradation in the seriousness of the conditions; a, study of the code-vectors profiles shows also a neat gradation of the negative modalities, which is an interesting result in that it could be usable as a basis for a "weighting" of the items.

In a second classification, the position of the modalities are represented on a 10 x 10 grid (Figure 2.a). At first glance, the map shows two regions, with the « neutral » modalities grouped at the top and on the right of the grid, and the « negative » ones at the bottom and on the left. This division is confirmed when we look at the representation of the distances between classes (Figure 2.b). In the bottom left cell we find again the characteristics of a poor standard of the dwelling, and generally, the bottom line corresponds to attributes of rather bad living conditions. The seriousness of the characteristics tends to decrease when we go towards the top of the map, but it is also interesting to notice that there are many modalities that have no immediate neighbor, suggesting that the hypothesis of cumulative difficulties could be somewhat reductive. If we group now the classes according to their closeness on the map, we obtain 3 groups of "negative" modalities: a first one grouping the characteristics of very serious living conditions (low standard dwelling, absence of very common durable goods, and privations in elementary consumptions), a second one corresponding mainly to other problems relating to the dwelling and environmental disadvantages, and a last group, at the "frontier" between "negative" and "neutral" conditions, characterized by the inability to afford one week's holidays or replacing worn out furniture that suggest a particular status for these items. These results are consistent with those obtained with a MCA (Figure 3).

## 4. Classes of households

We try now to classify the households according to their living conditions characteristics. First, in order to use the full qualitative resource of the initial information, we use in input the coordinates of the observations after a MCA (this corresponds to a transformation of the responses into quantitative values). Then we use only the partial scores calculated by domain of the living conditions, each household being then characterized by 5 scores.

### *4.1. Classification using the whole set of items*

The inputs are now the coordinates of the observations, resulting from a "traditional" MCA. The Kohonen algorithm is used to classify the observations in a 8 X 8 grid (since we have about 6000 observations, it could give about 100 observations by class). The 64 classes obtained are then grouped into 5 super-classes (SC1 to SC5), using a hierarchical classification. Figure 4 shows the Kohonen map obtained, and indicates the super-classes. For the analysis, another clustering of 3 super-classes - which results in keeping SC1 as one group A, adding SC2 and SC3 into another group B, and SC4 and SC5 into a third group C - is also used.

The analysis of the classes characteristics (table 1) shows some neat differences between group A and the two others: as for the characteristics of material living conditions only, group A appears to be in the best situation; this is also the case when we look at their monetary poverty rate, which is lower than on average, and their "subjective" living conditions, that appear easier than on average. But there are also significant differences between groups B and C: the households in group B are mainly disadvantaged in the domains of durables and privations, while those in group C are mostly concerned with low standard dwellings. The households in group B are also those facing the highest proportion of environmental problems, which is consistent with the proportion living in large structures and big cities; it is also the group where we find the highest proportion of lone parents households. If we look at their other characteristics, the monetary poverty rate is higher in group B than in group C, and "subjective" living conditions more often declared "very difficult". Going back to the five super-classes for more detail, we note that in fact, there are strong differences within groups B and C. Within group B, SC2 appears to suffer from very serious deprivations (in food and clothes), while group SC3 is mostly concerned by a lack of durables. Within group C, we have two sub-groups living in low standard dwellings, but there is a neat difference in the nature of missing items: mostly heating system and absence of a separate kitchen in SC4, and mostly basic comfort (bath or shower, hot running water and indoor toilets) in SC5; SC5 is also the class where we find the smallest average income and the highest proportion of poor. SC5, which counts for 2% of the households, is also characterized by a high proportion of households counting only one person, older than on average and living in rural areas.

These results tend to confirm that differentiations are not only in the level, but also in the nature of the disadvantages.

*4.2. Classification using the partial scores*

Our purpose now is to use a quantitative measurement of living conditions, but also to introduce some of the qualitative dimension of the information. We start from the simplest way to count (bad)living conditions "level", which consists in calculating a score of "bad points" for the whole set of items and setting one or several thresholds. The problem is then to determine where to put the thresholds to obtain consistent classes. This question is discussed in Lollivier & Verger (1997); they propose to use the rate of monetary poverty to set a "poverty line" for living conditions: the score that separates "good" and "bad" living conditions is the one for which the cumulative percent of the distribution of the score is equal (or close) to the monetary poverty rate. Applied to our data, the monetary poverty rate is 10,7% and at this cumulative percent of the score's distribution, the score equals 9 (table 2); so we'll say that living conditions are "good" while the score is under 9, and "bad" from a score of 9. One problem with this method is that it cuts the population into one group having "good" living conditions and one group having "bad" living conditions whatever the items; for example, the total value of living in a dwelling without hot running water and not being able to afford a meal of meat/fish/chicken every second day and not being able to afford buying new clothes is the same as having noisy neighbors and not having a VCR and not having a micro-wave oven. So we have tried in what follows to obtain classes defined using the qualitative dimension of the information and independently from an exogenous threshold.

For this, we have classified the observations according to their partial scores (one score by domain) of living conditions. We have used the Kohonen algorithm to classify the households along a one-dimension map of five classes (Figure 5). The "progression" appears neatly going from best to worse living conditions, but not only in terms of level: each of the 5 domains of living conditions appear to contribute more or less to the households positions (table 3). For example, class 3 is mostly characterized by environmental problems, this being also what makes the difference between classes 1 and 2. Then in classes 4 and 5, the partial scores are above average in almost all the domains of living conditions, except durables in class 4 and environmental problems in class 5.

If we look closer at the other characteristics of the households, we find that classes 4 and 5 have the highest proportions of persons living alone or lone parents, that classes 3 and 4 are more often living in large structures than on average, which is consistent with the proportion of environmental problems. More in the detail of the items, it appears clearly that the households of class 5 combine the highest proportions of difficulties in all the items of almost all the domains; the main difference with the households of class 4 being in the items relating to the standard of the dwelling. Class 3 is neatly disadvantaged on all the items relating to environmental quality, as for class 2, but with a smaller intensity and with almost no other material difficulty. Class 1 is clearly the one benefiting for the "best" material living conditions.

If we compare now the result from using the total score and the threshold set at 9 items and the result of our classification, a first main difference is in the proportion of households who will be said having bad living conditions: respectively 10,8% and (at least) 15%, if we consider in this case only the class 5, and up to almost 25% if we add classes 4 and 5. A second difference is that having more than two classes (one "good" and one "bad") allows to better differentiate between households who may have similar total score, but are different by the domains of their difficulties.

## 5. Conclusion

This paper was a first attempt to use the Kohonen algorithm to classify households' living conditions. It has proved a useful tool, in that it allows to construct classes that are neatly distinct and rather easily interpretable. In a first classification using the full qualitative information, we obtain classes of households which separates one class having "satisfying" material living conditions, and several classes of households who face more or less difficulties in almost all the domains considered; this classification makes appear neatly a small group of households characterized firstly by the lack of a basic comfort in their dwelling and the absence of some very common durable goods (telephone and television), as opposed to others groups characterized rather by serious deprivations in food and clothing. As for the results of the second classification based on the partial scores suggest strongly that poverty in terms of living conditions is not only a question of "level", but also a question of what contributes to the level.

The classifications show also that the most serious lacks or privations are often associated with a rather high score of (bad) living conditions, and a low income. Nevertheless, because they use the qualitative dimension of the information, they show some associations of characteristics that will define particular groups of households, that do not appear if we "measure" living conditions with a score.

**Acknowledgement**

The SAS programs used for Kohonen classifications are due to Patrick LETREMY (SAMOS/MATISSE, Université de Paris 1)

**References**

Cottrell, M. & Ibbou, S. (1995) : Multiple correspondence analysis of a crosstabulation matrix using the Kohonen algorithm, *Proc. ESANN'95*, M.Verleysen Ed., Editions D Facto, Bruxelles, 27-32.

Cottrell, M., Fort, J.C. & Pagès, G. (1998) : Theoretical aspects of the SOM Algorithm, *Neurocomputing,* 21, p. 119-138.

Cottrell, M. & Rousset, P. (1997) : The Kohonen algorithm : a powerful tool for analysing and representing multidimensional quantitative et qualitative data, *Proc. IWANN'97*, Lanzarote.

Dickes, P (1994), *Ressources financières, bien-être subjectif et conditions d'existence*, in F.Bouchayer (ed.) Trajectoires sociales et inégalités, Eres.

Kaski, S. (1997) : Data Exploration Using Self-Organizing Maps, *Acta Polytechnica Scandinavia*, 82.

Kohonen, T. (1984, 1993) : *Self-organization and Associative Memory*, 3$^{rd}$ Ed., Springer.

Kohonen, T. (1995) : Self-Organizing Maps, *Springer Series in Information Sciences* Vol 30, Springer.

Lollivier, S & Verger, D (1997), Pauvreté d'existence, monétaire ou subjective sont distinctes, *Economie & statistique* n°308/309/310.

Mack, J & Lansley, S (1984), *Poor Britain*, Allen & Unwin.

Mayer, S.E & Jencks, C (1989), *Poverty and the distribution of material hardship*, Journal of Human Resources vol.24.

Nolan, B & Whelan, C.T. (1996), *Resources, deprivation and poverty*, Clarendon Press.

Townsend, P (1979), *Poverty in the United-Kingdom*, Penguin.

Townsend, P (1989), Deprivation, *Journal of Social Policy*, vol.16.

**Figure 1 – Classification of the modalities in a one-dimension Kohonen map**

| CLB1 PHO1<br>CLW1 TEV1<br>CLE1 CAR1<br>PLT1 NCHF1<br>FMO1 NVIP1 | CLC1 EVB1<br>CLR1 VCR1<br>PLF1 LV1<br>PLM1 NAMI1<br>PLE1 NVET1<br>PLS1 | EB1<br>EP1<br>EV1<br>NMOB1<br>NVAC1 | EV0<br>NMOB0<br>NVAC0 | ALL THE OTHER « NEUTRAL » MODALITIES |
|---|---|---|---|---|
| 1 | 2 | 3 | 4 | 5 |

**Figure 2.a – Classification of the modalities in a 10 X 10 Kohonen map**

|  | EB1<br>EP1 |  | LV0 |  | FMO0<br>VCR0 |  | CAR0 |  | CLB0<br>CLW0<br>CLE0 |
|---|---|---|---|---|---|---|---|---|---|
| EVB1 |  | EV1 |  |  | NCHF0<br>NVIP0 |  | PLT0<br>PHO0<br>TEV0 |  |  |
|  |  |  | NAMI0 | NVET0 |  | EVB0 |  |  | CLC0<br>PLE0<br>PLS0 |
| CLB1<br>PLF1 |  | PLM1 |  |  | EB0<br>EP0 |  | CLR0 |  |  |
|  |  |  | NMOB0 | NVAC0 |  |  |  |  | PLF0<br>PLM0 |
| PLS1 |  | PLE1 |  |  | EV0 |  |  |  |  |
|  |  |  | VCR1 |  |  | NMOB1 |  | NVAC1 |  |
| PLT1 |  | CLC1 |  |  | LV1 |  |  |  |  |
|  |  |  |  | FMO1 |  |  |  | NAMI1<br>NVET1 |  |
| CLB1<br>CLW1<br>CLE1 |  | PHO1<br>TEV1 | CAR1 |  | NCHF1 |  | NVIP1 |  |  |

**Figure 2.b – Distances between the Kohonen classes of modalities**

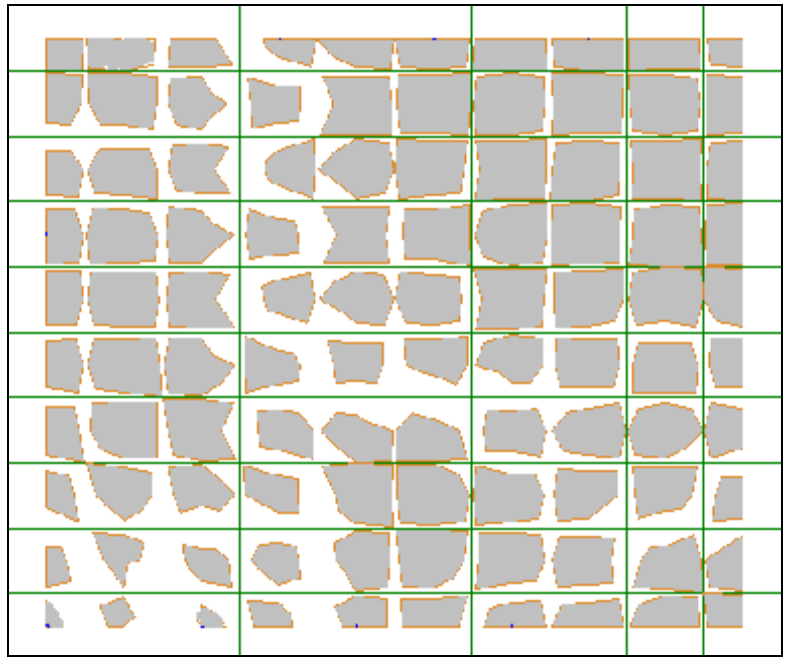

**Figure 3 – Multiple correspondence analysis of the modalities**
(projection of axes 1 – 2)

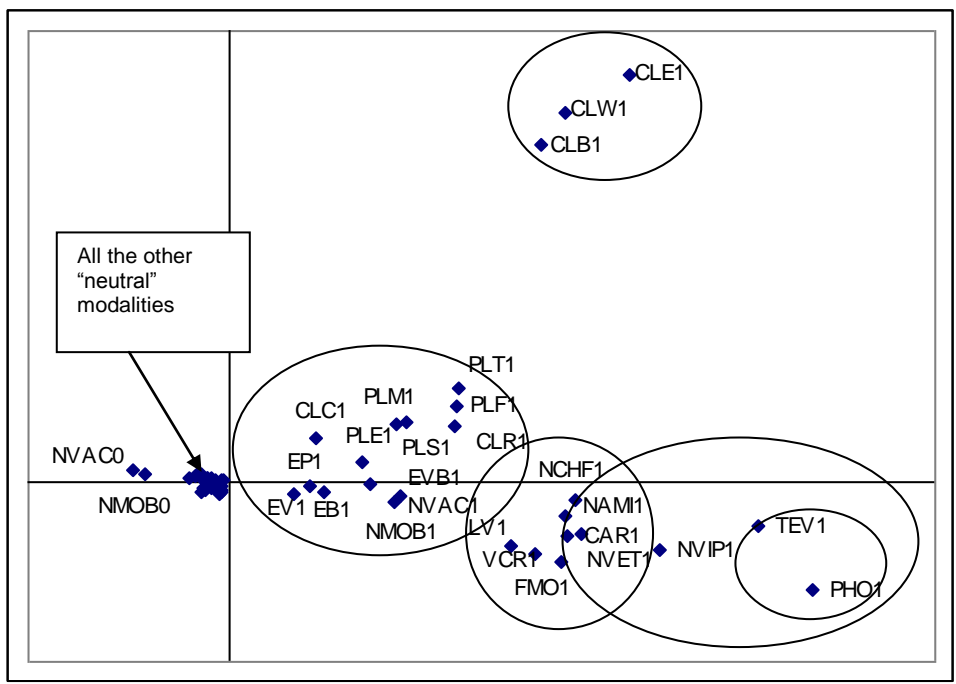

**Figure 4 – Classification of the households based on the full set of items: distances between the classes and representation of the super-classes**
*(coordinates after a MCA)*

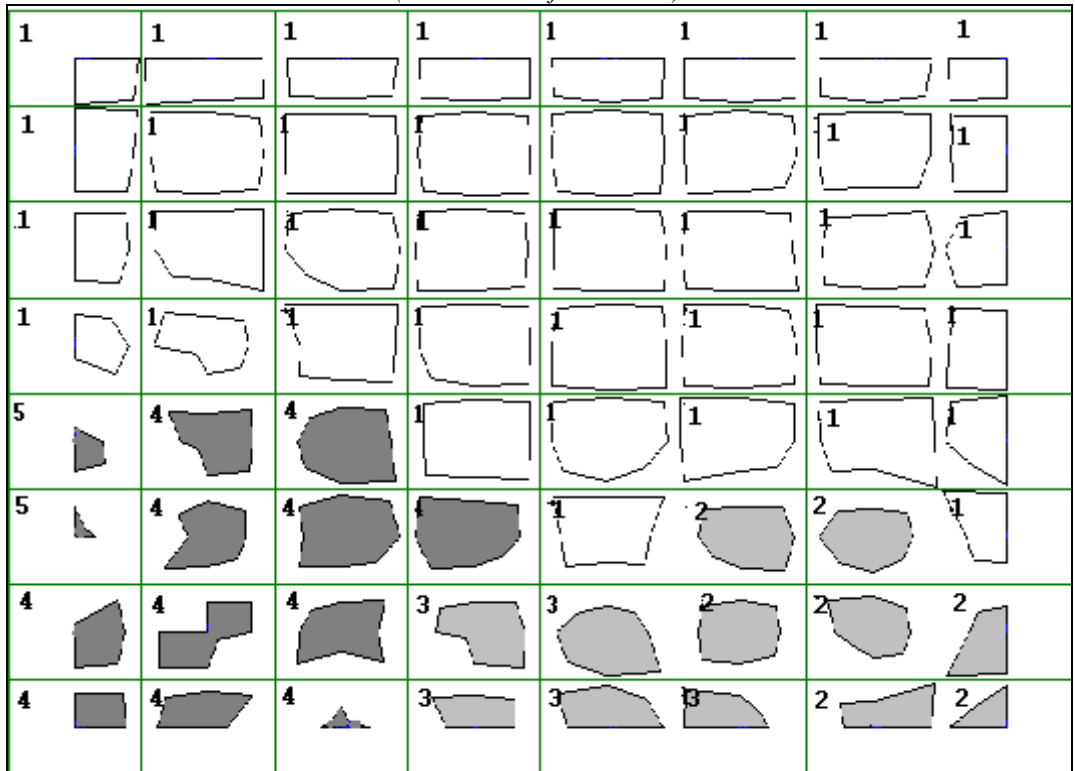

*Note: the numbers correspond to classes SC1 to SC5, and the shades of gray to groups A, B and C referred to in the text.*

**Table 1**
**a – Characteristics of the classes of households based on the full set of items**

|  | A | B | | | C | | | All |
|---|---|---|---|---|---|---|---|---|
|  | 1 | 2 | 3 | 2+3 | 4 | 5 | 4+5 |  |
| Proportion of each class | 70,6 | *8,7* | *5,5* | 14,2 | *13,3* | *1,9* | 15,2 | 100,0 |
| Average total score | 2,2 | **9,3** | **5,4** | 7,8 | **4,7** | **9,1** | **5,3** | 3,5 |
| Average partial scores: | | | | | | | | |
|   dwelling 1 | 0,1 | *0,3* | *0,3* | 0,3 | **1,0** | **2,8** | **1,2** | 0,3 |
|   dwelling 2 | 0,4 | **0,7** | **0,8** | 0,7 | **0,8** | **1,2** | **0,9** | 0,5 |
|   environment | 0,6 | **0,9** | **0,9** | 0,9 | *0,7* | *0,6* | *0,7* | 0,7 |
|   durables | 0,1 | **0,7** | **2,1** | **1,2** | *0,6* | *0,7* | *0,6* | 0,4 |
|   deprivations | 0,6 | **3,7** | **1,7** | **2,9** | *1,1* | **2,3** | *1,2* | 1,0 |
| Average income per C.U. | 8378 | **5103** | **5457** | 5240 | **6797** | **4565** | 6518 | 7650 |
| Proportion of poor | 6,1 | **27,6** | **19,8** | 24,6 | **17,0** | **33,3** | 19,0 | 10,7 |
| "Subjective" living conditions: | | | | | | | | |
|   with great difficulty | 2,0 | **25,7** | **14,1** | 21,2 | *7,2* | **13,9** | 8,0 | 5,7 |
|   with difficulty | 9,6 | **29,8** | **20,3** | 26,1 | *12,4* | *10,7* | 12,2 | 12,3 |
|   with some difficulty | 27,6 | **34,8** | **40,1** | 36,8 | **31,1** | **41,0** | 32,3 | 29,6 |
|   fairly easily | **44,6** | *8,6* | *22,0* | 13,8 | *38,1* | *27,1* | 36,8 | 39,0 |
|   easily and very easily | **16,1** | *1,3* | *3,4* | 2,1 | *11,2* | *7,4* | 10,7 | 13,3 |
| Type of dwelling | | | | | | | | |
| House, isolated | **43,5** | *28,5* | *27,4* | 28,1 | *33,2* | **41,5** | 34,3 | 39,9 |
| House, in a neighborhood | 21,5 | *20,3* | *22,0* | 21,0 | *18,8* | **36,6** | 21,0 | 21,4 |
| Structure <10 units | 11,8 | *15,7* | *15,8* | 15,7 | **19,2** | *13,8* | 18,5 | 13,4 |
| Structure >=10 units | 22,4 | **34,6** | **34,5** | 34,5 | **26,7** | *4,1* | 23,9 | 24,3 |
| Other | 0,8 | *0,9* | *0,3* | 0,7 | **2,1** | **4,1** | 2,3 | 1,0 |
| Type of location | | | | | | | | |
| Rural town | 27,5 | *24,4* | *26,8* | 25,4 | *28,8* | **45,5** | 30,9 | 27,7 |
| City <10000 inh | 11,2 | *8,7* | *10,7* | 9,5 | *11,4* | *6,5* | 10,8 | 10,9 |
| 10000 to <100000 inh | 19,0 | *20,9* | *21,2* | 21,0 | *20,9* | *17,1* | 20,4 | 19,5 |
| 100000 to <2000000 inh | 28,1 | **32,6** | **30,2** | 31,7 | *27,1* | *26,0* | 26,9 | 28,5 |
| Paris area | 14,2 | *13,4* | *11,0* | 12,5 | *11,9* | *4,9* | 11,0 | 13,5 |
| Type of household | | | | | | | | |
| one person household | 21,7 | **33,0** | *27,1* | 30,7 | **33,8** | **48,8** | 35,7 | 25,1 |
| couple without child | 28,2 | *17,3* | *21,8* | 19,0 | *25,4* | *15,5* | 24,2 | 26,3 |
| couple with child(ren) | 40,1 | *30,3* | *35,6* | 32,4 | *31,9* | *19,5* | 30,4 | 37,5 |
| lone parent family | 6,3 | **16,0** | **10,7** | 14,0 | *5,1* | *8,9* | 5,6 | 7,3 |
| other type | 3,8 | *3,4* | *4,8* | 3,9 | *3,7* | **7,3** | 4,2 | 3,8 |
| Total number of persons in the household | 2,6 | *2,5* | **2,7** | 2,6 | *2,4* | *2,0* | 2,3 | 2,6 |
| Age of the adults | 47,9 | *46,5* | *41,9* | 44,7 | *40,6* | **58,2** | 42,8 | 46,7 |
| Number of persons aged <17 | 0,6 | *0,6* | *0,7* | 0,6 | *0,6* | *0,3* | 0,5 | 0,6 |

*Note: the characteristics that are significantly over-represented are in bold*

**b – Detailed proportions for the full set of items**

|  | A | B | | | C | | | All |
|---|---|---|---|---|---|---|---|---|
|  | 1 | 2 | 3 | 2+3 | 4 | 5 | 4+5 |  |
| no bath or shower | 2,7 | 0,0 | 0,0 | 0,0 | 0,7 | **71,5** | 9,6 | 3,4 |
| no hot running water | 0,0 | 0,0 | 0,0 | 0,0 | 0,2 | **100,0** | 12,7 | 1,9 |
| not indoor toilet | 1,8 | 0,9 | 0,6 | 0,8 | 2,3 | **49,6** | 8,2 | 2,6 |
| inappropriate heating system | 1,9 | 17,5 | 16,7 | 17,2 | **45,5** | 30,1 | 43,6 | 9,5 |
| no separate kitchen | 1,0 | 10,9 | 8,2 | 9,8 | **52,0** | 26,0 | 48,8 | 10,4 |
| rot in window frames | 6,6 | 15,3 | 13,6 | 14,6 | 16,8 | **30,1** | 18,5 | 9,6 |
| damp walls, floor, etc. | 11,5 | 19,8 | 21,2 | 20,3 | 22,2 | **31,7** | 23,4 | 14,6 |
| leaky roof | 3,6 | 6,4 | 6,8 | 6,6 | 7,9 | **17,9** | 9,2 | 4,8 |
| shortage of space | 9,2 | 19,1 | **20,6** | 19,7 | 18,6 | 15,5 | 18,2 | 12,0 |
| too dark, not enough light | 6,7 | 13,6 | 13,0 | 13,3 | 16,3 | 21,1 | 16,9 | 9,2 |
| noise from outside | 16,0 | **26,4** | **26,3** | 26,3 | 23,3 | 17,1 | 22,6 | 18,5 |
| noise from neighbors | 8,3 | **13,2** | **15,5** | 14,1 | **14,2** | 9,8 | 13,6 | 9,9 |
| pollution | 13,7 | **19,3** | **19,5** | 19,3 | 14,9 | **19,5** | 15,5 | 14,8 |
| vandalism, lack of safety | 20,7 | **31,4** | **27,1** | 29,7 | 21,5 | 16,3 | 20,8 | 22,0 |
| no telephone | 0,0 | 0,0 | 0,0 | 0,0 | **9,1** | **4,1** | 8,4 | 1,3 |
| no color TV | 0,0 | 0,0 | 0,6 | 0,2 | **8,5** | **4,9** | 8,0 | 1,3 |
| no car | 3,1 | 14,6 | 18,9 | 16,3 | 10,0 | 10,6 | 10,1 | 6,0 |
| no microwave | 0,2 | 10,3 | **100,0** | 45,0 | 7,2 | 14,6 | 8,1 | 7,8 |
| no VCR | 4,1 | 20,7 | **39,8** | 28,1 | 15,6 | 17,1 | 15,8 | 9,3 |
| no dishwasher | 5,4 | 20,5 | **48,9** | 31,5 | 11,6 | 15,5 | 12,1 | 10,1 |
| Cannot afford |  |  |  |  |  |  |  |  |
| meat/fish/chicken every second day if wanted | 0,2 | **46,9** | 0,0 | 28,7 | 3,6 | 16,3 | 5,2 | 5,0 |
| new clothes | 0,4 | **76,1** | 15,0 | 52,5 | 6,6 | 25,2 | 8,9 | 9,1 |
| keeping the dwelling warm enough | 3,1 | 28,0 | 11,9 | 21,8 | 9,9 | 20,3 | 11,2 | 6,9 |
| replacing worn out furniture | 26,2 | *92,0* | *68,6* | *83,0* | *43,9* | *61,0* | *46,0* | 37,3 |
| a week's holiday away once a year | 24,5 | *79,0* | *57,9* | *70,8* | *36,7* | *69,1* | *40,8* | 33,5 |
| having friends/family for a drink or meal at least once a month | 5,8 | **48,5** | 20,1 | 37,5 | 7,3 | 33,3 | 10,6 | 11,0 |

**Table 2 – Setting a poverty line of living conditions**

| SCORE | Percent | Cumulative Percent | |
|---|---|---|---|
| 23 | 0,0 | 100, | 0,0 |
| 22 | 0,0 | 100 | 0,0 |
| 21 | 0,0 | 100 | 0,0 |
| 20 | 0,0 | 100 | 0,0 |
| 19 | 0,1 | 99,9 | 0,1 |
| 18 | 0,1 | 99,9 | 0,2 |
| 17 | 0,1 | 99,8 | 0,3 |
| 16 | 0,2 | 99,6 | 0,5 |
| 15 | 0,6 | 99,4 | 1,1 |
| 14 | 0,6 | 98,8 | 1,7 |
| 13 | 0,9 | 98,2 | 2,6 |
| 12 | 1,4 | 97,3 | 4,0 |
| 11 | 1,6 | 95,9 | 5,6 |
| 10 | 2,3 | 94,4 | 7,9 |
| **9** | 2,9 | 92,1 | **10,8** |
| 8 | 3,5 | 89,2 | 14,3 |
| 7 | 3,9 | 85,7 | 18,2 |
| 6 | 5,7 | 81,8 | 23,9 |
| 5 | 6,5 | 76,2 | 30,4 |
| 4 | 9,0 | 69,6 | 39,4 |
| 3 | 8,6 | 60,6 | 48,0 |
| 2 | 14,5 | 52,0 | 62,5 |
| 1 | 14,3 | 37,5 | 76,8 |
| 0 | 23,2 | 23,2 | 100 |

**Figure 5 – Means of 5 Kohonen classes of households based on the partial scores**

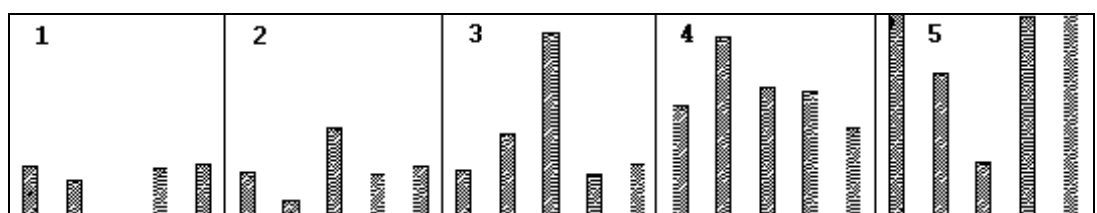

*Note: the graph shows the standardized values of the means of the partial scores; the order is as follows: dwelling1, dwelling2, environment, privations, durables.*

**Table 3 - 5 Kohonen classes of households based on partial scores**

|  | Basic score | | Kohonen classes | | | | | All |
|---|---|---|---|---|---|---|---|---|
|  | 0 | 1 | 1 | 2 | 3 | 4 | 5 |  |
| Proportion of each class | 89,2 | 10,8 | 48,7 | 12,5 | 14,0 | 9,9 | 15,0 | 100,0 |
| Average partial scores: | | | | | | | | |
|    dwelling 1 | 0,2 | 0,9 | 0,1 | 0,1 | 0,1 | **0,5** | **1,0** | 0,3 |
|    dwelling 2 | 0,4 | 1,4 | 0,2 | 0,0 | 0,6 | **1,5** | **1,2** | 0,5 |
|    environment | 0,6 | 1,4 | 0,0 | 1,0 | **2,1** | **1,5** | 0,6 | 0,7 |
|    durables | 0,3 | 1,1 | 0,1 | 0,1 | 0,1 | 0,5 | **1,6** | 0,4 |
|    deprivations | 0,7 | 3,8 | 0,6 | 0,5 | 0,5 | **1,8** | **3,0** | 1,0 |
| Average total score | 2,5 | 11,3 | 1,5 | 2,0 | 3,7 | **7,0** | **8,8** | 3,5 |
| Average income per C.U. | 8021 | 4591 | 8251 | 9092 | 8671 | 6026 | 4613 | 7650 |
| Proportion of poor | 7,8 | 34,7 | 6,3 | 3,9 | 4,4 | **15,3** | **33,4** | 10,7 |
| "Subjective" living conditions: | | | | | | | | |
|    with great difficulty | 2,6 | 30,9 | 1,5 | 1,3 | 2,7 | **9,4** | **23,4** | 5,7 |
|    with difficulty | 10,5 | 28,0 | 8,1 | 7,9 | 9,5 | **20,9** | **26,7** | 12,3 |
|    with some difficulty | 29,3 | 32,7 | 27,3 | 25,8 | 28,6 | **39,5** | **34,9** | 29,6 |
|    fairly easily | 42,8 | 7,5 | **46,4** | **46,8** | **43,8** | 25,5 | 13,1 | 39,0 |
|    easily and very easily | 14,8 | 1,0 | **16,7** | **18,3** | **15,4** | 4,7 | 1,9 | 13,3 |
| Details | | | | | | | | |
| no bath or shower | 2,1 | 13,8 | 1,4 | 0,3 | 0,2 | 2,3 | **15,9** | 3,4 |
| no hot running water | 1,0 | 9,9 | 0,4 | 0,5 | 0,2 | 1,7 | **10,0** | 1,9 |
| not indoor toilet | 1,6 | 11,2 | 0,7 | 0,4 | 0,3 | 3,3 | **12,4** | 2,6 |
| no separate kitchen | 7,3 | 36,5 | 3,0 | 4,4 | 5,0 | 24,2 | **35,6** | 9,5 |
| inappropriate heating system | 8,5 | 18,2 | 7,0 | 3,9 | 4,1 | 14,1 | **24,6** | 10,4 |
| rot in window frames | 6,9 | 31,5 | 2,1 | 0,0 | 12,1 | **31,9** | 24,6 | 9,6 |
| damp walls, floor, etc. | 11,5 | 40,4 | 5,9 | 0,0 | 17,6 | **42,3** | 33,8 | 14,6 |
| leaky roof | 3,4 | 17,1 | 1,6 | 0,0 | 3,9 | **16,1** | 12,8 | 4,8 |
| shortage of space | 10,1 | 28,2 | 5,5 | 0,0 | 17,9 | **33,0** | 23,8 | 12,0 |
| too dark, not enough light | 7,0 | 27,4 | 3,3 | 0,0 | 10,4 | **27,7** | 22,7 | 9,2 |
| noise from outside | 15,8 | 40,4 | 0,0 | 21,2 | **63,6** | 45,5 | 16,2 | 18,5 |
| noise from neighbors | 7,8 | 27,4 | 0,0 | 10,4 | **31,5** | 24,2 | 12,2 | 9,9 |
| pollution | 13,0 | 29,7 | 0,0 | 18,1 | **53,7** | 31,4 | 12,6 | 14,8 |
| vandalism, lack of safety | 19,3 | 43,8 | 0,0 | 50,2 | **58,6** | 45,9 | 19,8 | 22,0 |
| no telephone | 0,7 | 6,2 | 0,1 | 0,1 | 0,0 | 1,4 | **7,4** | 1,3 |
| no color TV | 0,9 | 4,9 | 0,1 | 0,0 | 0,1 | 0,8 | **7,6** | 1,3 |
| no car | 4,2 | 20,8 | 1,4 | 1,7 | 2,0 | 9,1 | **26,3** | 6,0 |
| no microwave | 5,8 | 24,4 | 2,3 | 1,6 | 1,7 | 9,7 | **35,3** | 7,8 |
| no VCR | 7,0 | 27,9 | 3,4 | 2,0 | 3,2 | 9,7 | **39,9** | 9,3 |
| no dishwasher | 7,7 | 29,7 | 3,8 | 3,2 | 3,7 | 14,4 | **39,8** | 10,1 |
| Cannot afford | | | | | | | | |
| meat/fish/chicken every second day if wanted | 1,4 | 35,2 | 1,0 | 1,0 | 0,8 | 8,0 | **23,4** | 5,0 |
| new clothes | 3,7 | 53,6 | 2,9 | 2,0 | 2,9 | 17,7 | **35,5** | 9,1 |
| keeping the dwelling warm enough | 3,0 | 39,1 | 1,6 | 2,4 | 2,1 | 12,8 | **28,7** | 6,9 |
| replacing worn out furniture | 30,3 | 95,3 | *25,0* | *20,1* | *22,0* | *65,0* | *87,5* | 37,3 |
| a week's holiday away once a year | 4,8 | 61,8 | *23,5* | *17,3* | *15,7* | *54,4* | *82,5* | 33,5 |
| having friends/family for a drink or meal at least once a month | 26,5 | 92,0 | 3,7 | 3,0 | 3,2 | 20,9 | **42,0** | 11,0 |

**APPENDIX Table A – Detailed frequencies of the items**

| Domain / Items | Answer | Variable | code | % |
|---|---|---|---|---|
| **Basic comfort of the dwelling** | | | | |
| bath or shower | yes | CLB | 0 | 96.6 |
| | no | | 1 | 3.4 |
| hot running water | yes | CLE | 0 | 98.1 |
| | no | | 1 | 1.9 |
| indoor flushing toilet | yes | CLW | 0 | 97.4 |
| | no | | 1 | 2.6 |
| appropriate heating | yes | CLR | 0 | 89.6 |
| | no | | 1 | 10.4 |
| separate kitchen | yes | CLC | 0 | 90.5 |
| | no | | 1 | 9.5 |
| **Problems in the dwelling** | | | | |
| rot in windows frames | no | PLF | 0 | 90.4 |
| | yes | | 1 | 9.6 |
| damp walls, floors… | no | PLM | 0 | 85.4 |
| | yes | | 1 | 14.6 |
| leaky roof | no | PLT | 0 | 95.2 |
| | yes | | 1 | 4.8 |
| shortage of space | no | PLE | 0 | 88.0 |
| | yes | | 1 | 12.0 |
| too dark, not enough light | no | PLS | 0 | 90.8 |
| | yes | | 1 | 9.2 |
| **Environment** | | | | |
| noise from outside | no | EB | 0 | 81.5 |
| | yes | | 1 | 18.5 |
| noise from neighbors | no | EVB | 0 | 90.1 |
| | yes | | 1 | 9.9 |
| pollution from traffic or industry | no | EP | 0 | 85.2 |
| | yes | | 1 | 14.8 |
| vandalism or crime in the area | no | EV | 0 | 78.0 |
| | yes | | 1 | 22.0 |
| **Durables** | | | | |
| telephone | yes | PHO | 0 | 98.7 |
| | no | | 1 | 1.3 |
| color TV | yes | TEV | 0 | 98.7 |
| | no | | 1 | 1.3 |
| car | yes | CAR | 0 | 94.0 |
| | no | | 1 | 6.0 |
| micro-wave oven | yes | FMO | 0 | 92.2 |
| | no | | 1 | 7.8 |
| VCR | yes | VCR | 0 | 90.7 |
| | no | | 1 | 9.3 |
| dish-washer | yes | LV | 0 | 89.9 |
| | no | | 1 | 10.1 |
| **Deprivations** | | | | |
| buying meat-chicken-fish every second day | yes | NVIP | 0 | 95.0 |
| | no | | 1 | 5.0 |
| buying new clothes | yes | NVET | 0 | 90.9 |
| | no | | 1 | 9.1 |
| keeping the dwelling adequately warm | yes | NCHF | 0 | 93.1 |
| | no | | 1 | 6.9 |
| replacing worn out furniture | yes | NMOB | 0 | 62.7 |
| | no | | 1 | 37.3 |
| paying for a week's holiday away once a year | yes | NVAC | 0 | 66.5 |
| | no | | 1 | 33.5 |
| having friends or family for a meal once a month | yes | NAMI | 0 | 89.0 |
| | no | | 1 | 11.0 |

**APPENDIX Table B – General descriptors**

|  |  | Variable |  | Frequencies (%) | Means |
|---|---|---|---|---|---|
| Type of dwelling | House, isolated | LOGT | 1 | 39.9 |  |
|  | House, in a neighborhood |  | 2 | 21.4 |  |
|  | Structure <10 units |  | 3 | 13.4 |  |
|  | Structure >=10 units |  | 4 | 24.3 |  |
|  | Other |  | 5 | 1.0 |  |
| Location | Rural town | TUR | 0 | 27.7 |  |
|  | City <10000 inh |  | 1 | 10.9 |  |
|  | 10000 to <100000 inh |  | 2 | 19.5 |  |
|  | 100000 to <2000000 inh |  | 3 | 28.5 |  |
|  | Paris area |  | 4 | 13.5 |  |
| Type of household (children taken into account if <25 years old) | one person household | TYM | 0 | 25.1 |  |
|  | couple without child |  | 1 | 26.3 |  |
|  | couple with child(ren) |  | 2 | 37.5 |  |
|  | lone parent family |  | 3 | 7.3 |  |
|  | other type |  | 4 | 3.8 |  |
| Total number of individuals |  | NBTOT |  |  | 2.6 |
| Number of children <17 years old |  | NB17 |  |  | 0.6 |
| Mean age of persons at least 17 years old |  | AGEM |  |  | 46.7 |
| "Subjective" living conditions | very difficult | SLS | 1 | 5.7 |  |
|  | difficult |  | 2 | 12.3 |  |
|  | rather difficult |  | 3 | 29.6 |  |
|  | rather comfortable |  | 4 | 39.0 |  |
|  | comfortable and very c. |  | 5 | 13.3 |  |
| Monthly income (FRF) per C.U.*(a)* |  | REVUC |  |  | 7650 |
| Monetary poverty | poor*(b)* | POOR | 1 | 10.7 |  |
|  | non poor |  | 2 | 89.3 |  |
| Score for the material living conditions | Total (all domains) | TOT4 |  |  | 3.5 |
|  | Dwelling, convenience | CLOGT |  |  | 0.3 |
|  | Dwelling, problems | PLOGT |  |  | 0.5 |
|  | Environment | ENVIR |  |  | 0.7 |
|  | Durables | DURAB |  |  | 0.4 |
|  | Deprivations | PRIV |  |  | 1.0 |

Source : Insee, European community households' panel, wave 3 (1996)
*(a)* Consumption Unit, using the following equivalence scale : 1 - 0.5 - 0.3
*(b)* Poverty threshold at 50 % of the median income per C.U.